\documentclass[12pt,reqno]{article}
\usepackage{graphicx}
\usepackage{amsmath}
\usepackage{amsthm}
\usepackage{latexsym,bm}
\usepackage{amssymb}
\usepackage{indentfirst}
\newtheorem{thm}{Theorem}[section]

\newtheorem{prop}[thm]{Proposition}
\newtheorem{rem}[thm]{Remark}

\numberwithin{equation}{section}

\begin{document}

\title{{\bf Some remarks on the Pigola-Rigoli-Setti version of the Omori-Yau maximum principle}}

\author{{\bf F. Fontenele\thanks{Work partially supported by CNPq
(Brazil)}\, and A.P. Barreto}}


\date{}
\maketitle

\begin{quote}
\small {\bf Abstract}. We prove that the  hypotheses in the
version of the Omori-Yau maximum principle that was given by
Pigola-Rigoli-Setti are logically equivalent to the assumption
that the manifold carries a $C^2$ proper function whose gradient
and Hessian (Laplacian) are bounded. In particular, this result
extends the scope of the original Omori-Yau principle, formulated
in terms of lower bounds for curvature.
\end{quote}

\begin{quote}

\end{quote}




\section{Introduction}

The celebrated Omori-Yau maximum principle (\cite{CY}, \cite{O},
\cite{Y1}) states that if $(M,g)$ is a complete Riemannian
manifold with sectional curvature (resp. Ricci curvature) bounded
below then, for every $f\in C^2(M)$ that is bounded above, there
exists a sequence $(q_k)$ in $M$ such that
\begin{eqnarray}\label{OY1}
f(q_k)>\sup_Mf-\frac{1}{k}\,,\;\;|\nabla
f|(q_k)<\frac{1}{k}\,,\;\;\textnormal{Hess}f(q_k)<\frac{1}{k}\,g(q_k)
\;\Big(\text{resp.}\;\Delta f(q_k)<\frac{1}{k}\Big),
\end{eqnarray}
for all $k\in\mathbb N$, where the third inequality above is in
the sense of quadratic forms.

\vskip5pt

Pigola-Rigoli-Setti (\cite{PRS}, Theorem 1.9) obtained a version
of the Omori-Yau maximum principle where the hypothesis that the
curvature is bounded below is replaced by the assumption that the
manifold admits a smooth function with special properties. More
precisely, they proved the following result:

\begin{thm}\label{OY2}
Let $(M,g)$ be a Riemannian manifold. Assume that there exist a
$C^2$ function $\gamma:M\to [0,+\infty)$, a compact set $K\subset
M$ and constants $A,B>0$ such that

\vskip5pt

(i) $\gamma$ is proper, i.e., $\gamma(x)\to +\infty$ \;as\;\;$x\to
\infty$,

\vskip3pt

(ii) $|\nabla\gamma|\leq
A\sqrt{\gamma}$\;\;\;\textnormal{on}\;\;$M\backslash K$,

\vskip3pt

(iii) $\textnormal{Hess}\,\gamma\leq B\sqrt{\gamma
G(\sqrt{\gamma}\,)}\,g$\;\;\Big(\textnormal{resp}.
$\Delta\gamma\leq B\sqrt{\gamma
G(\sqrt{\gamma}\,)}$\;\Big)\;\;\;\textnormal{on}\;\;$M\backslash
K$,

\vskip5pt

\noindent{where} $G:[0,+\infty)\to [0,+\infty)$ is a smooth
function satisfying
\begin{eqnarray}\label{OY2b}
G(0)>0,\;\;\;\;G'(t)\geq 0\;\textnormal{for all}\;t\geq 0,
\end{eqnarray}
\begin{eqnarray}\label{OY2a}
\int_0^{+\infty}\frac{dt}{\sqrt{G(t)}}=+\infty\;\;\;\;\textnormal{and}\;\;\;\;
\limsup_{t\to +\infty}\frac{tG(\sqrt{t})}{G(t)}<+\infty.
\end{eqnarray}
Then, for every $f\in C^2(M)$ that is bounded above, there exists
a sequence $(q_k)$ in $M$ satisfying (\ref{OY1}).
\end{thm}

The function theoretic approach to the Omori-Yau maximum principle
provided by Theorem \ref{OY2} has been applied by several 
authors to obtain results in different contexts (\cite{ABD},
\cite{AR}, \cite{PRS}). In Section 2 of this work, we show that
Theorem \ref{OY2} is logically equivalent to a more conceptual
statement:

\begin{thm}\label{OY3}
If a Riemannian manifold $(M,g)$ admits a $C^2$ function
$\phi:M\to\mathbb R$ satisfying, for some constants $C,D>0$,

\vskip5pt

(i) $\phi$ is proper,

\vskip3pt

(ii) $|\nabla\phi|\leq C$,

\vskip3pt

(iii) $\textnormal{Hess}\,\phi\leq Dg$ \;\;
\big(\textnormal{resp}. \;$\Delta\phi\leq D\big)$,

\vskip5pt

\noindent{then,} for every $f\in C^2(M)$ that is bounded above,
there exists a sequence $(q_k)$ in $M$ satisfying (\ref{OY1}).
\end{thm}


Under the assumption that a Riemannian manifold $M$ is complete
and has sectional curvature (resp. Ricci curvature) bounded below,
Schoen-Yau (\cite{SY}, Theorem 4.2) proved that, for every $p\in
M$, there exists a smooth function $\phi:M\to\mathbb R$
satisfying, for some constants $C,D>0$,
\begin{eqnarray}\label{OY5}
\phi(x)\geq d(x,p),\;\;|\nabla\phi|(x)\leq
C,\;\;\textnormal{Hess}\,\phi(x)(v,v)\leq D|\,v|^2
\;\;\;\;(\text{resp}. \;\Delta\phi(x)\leq D),\nonumber
\end{eqnarray}
for all $x\in M$ and $v\in T_xM$. In particular, $\phi$ satisfies
(i), (ii) and (iii) in the statement of Theorem \ref{OY3}. This
shows that Theorem \ref{OY3} (and so Theorem \ref{OY2}) is a
generalization of the original Omori-Yau maximum principle.

\section{The arguments.}

The proof below was inspired  by  the proof of a conceptual
refinement of the Omori-Yau maximum principle stated in \cite{FX}.

\vskip10pt

\noindent{\bf Proof of the Theorem \ref{OY3}.} We will give the
proof of the half of the theorem that refers to the Hessian. The
proof of the other half is entirely analogous and will be left to
the reader.

Multiplying $\phi$ by a positive constant if necessary, we can
assume that $C=D=1$. Let $f:M\to\mathbb R$ be a $C^2$ function
satisfying $\sup_Mf<+\infty$ and $(p_k)$ a sequence in $M$ such
that
\begin{eqnarray}\label{OY5a}
f(p_k)>\sup_Mf-\frac{1}{2k},\;\;\;k\in\mathbb N.
\end{eqnarray}
For each $k\in\mathbb N$, define a function $f_k:M\to\mathbb R$ by
\begin{eqnarray}\label{OY6}
f_k(x)=f(x)-\varepsilon_k [\phi(x)-\phi(p_k)],
\end{eqnarray}
where $\varepsilon_k=\min\{\eta_k,1/2k\}$ and
\begin{eqnarray}\label{OY7}
\eta_k=\begin{cases}1/2k[\phi(p_k)-\inf_M\phi],\;\;\;\text{if}\;\;\;\phi(p_k)>\inf_M\phi,\\
1/2k,\;\;\;\text{if}\;\;\;\phi(p_k)=\inf_M\phi.\end{cases}\nonumber
\end{eqnarray}
Since $\varepsilon_k>0$ and $f$ is bounded above, from (i) one
obtains that $f_k(x)\to -\infty$ as $x\to\infty$, and so $f_k$
attains a global maximum at some point $q_k\in M$. Hence, by
(\ref{OY6}) and the ordinary maximum principle,
\begin{eqnarray}\label{OY8}
0=\nabla f_k(q_k)=\nabla f(q_k)-\varepsilon_k\nabla\phi(q_k)
\end{eqnarray}
and
\begin{eqnarray}\label{OY9}
0\geq
\text{Hess}(f_k)(q_k)(v,v)=\text{Hess}f(q_k)(v,v)-\varepsilon_k\text{Hess}\,\phi(q_k)(v,v),\;\;\forall
v\in T_{q_k}M.
\end{eqnarray}
From (ii), (\ref{OY8}) and definition of $\varepsilon_k$, one
obtains
\begin{eqnarray}\label{OY10}
|\nabla f(q_k)|=\varepsilon_k |\nabla\phi(q_k)|\leq
\varepsilon_k\leq \frac{1}{2k}<\frac{1}{k}\,\cdot\nonumber
\end{eqnarray}
By (iii) and (\ref{OY9}), we have, for all $v\in T_{q_k}M$ with
$v\neq 0$,
\begin{eqnarray}\label{OY11}
\text{Hess}f(q_k)(v,v)\leq
\varepsilon_k\text{Hess}\,\phi(q_k)(v,v)\leq\varepsilon_k
|\,v|^2\leq\frac{1}{2k}|\,v|^2<\frac{1}{k}|\,v|^2.\nonumber
\end{eqnarray}
Since $f_k(p_k)=f(p_k)$, we also have
\begin{eqnarray}\label{OY12}
f(p_k)&=&f_k(p_k)\leq f_k(q_k)=f(q_k)-\varepsilon_k
[\phi(q_k)-\phi(p_k)]\nonumber\\&=&f(q_k)-\varepsilon_k
[\phi(q_k)-\inf_M\phi]-\varepsilon_k
[\inf_M\phi-\phi(p_k)]\nonumber\\&\leq& f(q_k)-\varepsilon_k
[\inf_M\phi-\phi(p_k)]\leq f(q_k)+\frac{1}{2k}\,\cdot\nonumber
\end{eqnarray}
Therefore, by (\ref{OY5a}),
\begin{eqnarray}\label{OY13}
f(q_k)\geq
f(p_k)-\frac{1}{2k}>\sup_Mf-\frac{1}{2k}-\frac{1}{2k}=\sup_Mf-\frac{1}{k},\nonumber
\end{eqnarray}
which completes the proof of the theorem.\qed

\vskip10pt


The fact that Theorem \ref{OY2} is equivalent to Theorem \ref{OY3}
is an immediate consequence of the following proposition.

\begin{prop}\label{OY13a}
A Riemannian manifold $(M,g)$ admits a function
$\gamma:M\to\mathbb R$ as in the statement of Theorem \ref{OY2} if
and only if it admits a function $\phi:M\to\mathbb R$ as in the
statement of Theorem \ref{OY3}.
\end{prop}

\noindent{\bf Proof.} Let $\gamma:M\to\mathbb R$ be a function as
in the statement of Theorem \ref{OY2}. Define a (smooth) function
$u:(0,+\infty)\to (0,+\infty)$ by
\begin{eqnarray}\label{OY14}
u(t)=\sqrt{tG(\sqrt{t})}.\nonumber
\end{eqnarray}
From (\ref{OY2b}), one obtains
\begin{eqnarray}\label{OY15}
u'(t)=\frac{G(\sqrt{t})+\frac{1}{2}\sqrt{t}G'(\sqrt{t})}{2\sqrt{tG(\sqrt{t})}}>0,\;\;\;
t>0.\nonumber
\end{eqnarray}
Given $C>\limsup_{t\to +\infty}tG(\sqrt{t})/G(t)$, there exists
$t_o>1$ such that
\begin{eqnarray}\label{OY16}
\frac{tG(\sqrt{t})}{G(t)}< C,\;\;\;t\geq t_o.\nonumber
\end{eqnarray}
From the above inequality and the fact that $G$ is non-decreasing,
we obtain
\begin{eqnarray}\label{OY17}
0<tG(0)\leq tG(\sqrt{t})< CG(t),\;\;\;t\geq t_o,\nonumber
\end{eqnarray}
and so
\begin{eqnarray}\label{OY18}
\frac{1}{\sqrt{tG(\sqrt{t})}}>\frac{1}{\sqrt{C}\sqrt{G(t)}},\;\;\;t\geq
t_o.\nonumber
\end{eqnarray}
Hence
\begin{eqnarray}\label{OY19}
\int_1^{+\infty}\frac{1}{u(s)}\,ds&=&\int_1^{+\infty}\frac{1}{\sqrt{sG(\sqrt{s})}}\,ds\geq
\int_{t_o}^{+\infty}\frac{1}{\sqrt{sG(\sqrt{s})}}\,ds\nonumber\\&\geq&\frac{1}{\sqrt{C}}
\int_{t_o}^{+\infty}\frac{1}{\sqrt{G(s)}}\,ds=+\infty,
\end{eqnarray}
where in the last equality we used (\ref{OY2a}).

Since $u(t)$ and $\sqrt{t}$ are non-decreasing, we can assume,
adding a positive constant if necessary, that $\gamma>0$ and that
(ii) and (iii) in the statement of Theorem \ref{OY2} holds on all
of $M$. Therefore, for all $x\in M$ and all $v\in T_xM$ we have
\begin{eqnarray}\label{OY20}
|\nabla\gamma(x)|\leq
A\sqrt{\gamma(x)}=\frac{A}{\sqrt{G(0)}}\sqrt{\gamma(x)G(0)}&\leq&\frac{A}{\sqrt{G(0)}}\sqrt{\gamma(x)G(\sqrt{\gamma(x)})}
\nonumber\\&\leq&\frac{A}{\sqrt{G(0)}}\,u(\gamma(x))
\end{eqnarray}
and
\begin{eqnarray}\label{OY21}
\text{Hess}\,\gamma(x)(v,v)\leq
B\sqrt{\gamma(x)G\big(\sqrt{\gamma(x)}\,\big)}\,|\,v|^2=Bu(\gamma(x))|\,v|^2.
\end{eqnarray}

Let $h:(0,+\infty)\to\mathbb R$ be defined by
\begin{eqnarray}\label{OY22}
h(t)=\int_1^t\frac{1}{u(s)}\,ds\nonumber
\end{eqnarray}
Since $u>0$ and $u'>0$, we have, for all $t>0$,
\begin{eqnarray}\label{OY23}
h'(t)=\frac{1}{u(t)}>0,\;\;\;\;h''(t)=-\frac{u'(t)}{u^2(t)}<0
.\end{eqnarray}

Let $\phi=h\circ\gamma$, so that
\begin{eqnarray}\label{OY23a}
\phi(x)=h(\gamma(x))=\int_1^{\gamma(x)}\frac{1}{u(s)}\,ds,\;\;\;x\in
M.
\end{eqnarray}
From (\ref{OY20}) and (\ref{OY23}), we have
\begin{eqnarray}\label{OY24}
|\nabla\phi(x)|=|\,h'(\gamma(x))\nabla\gamma(x)|=\frac{1}{u(\gamma(x))}|\nabla\gamma(x)|\leq
\frac{A}{\sqrt{G(0)}},\;\;\;\forall x\in M.
\end{eqnarray}
Using (\ref{OY21}) and (\ref{OY23}), we obtain, for all $x\in M$
and all $v\in T_xM$,
\begin{eqnarray}\label{OY25}
\text{Hess}\,\phi(x)(v,v)&=&h''(\gamma(x))\langle\nabla\gamma(x),v\rangle^2+h'(\gamma(x))\text{Hess}\gamma(x)(v,v)
\nonumber\\&\leq&\frac{1}{u(\gamma(x))}\text{Hess}\gamma(x)(v,v)\leq
B|\,v|^2.
\end{eqnarray}
Moreover, from (\ref{OY19}), (\ref{OY23a}) and the properness of
$\gamma$, we obtain that $\phi$ is proper. This concludes the
proof of the ``only if" part of the proposition. The ``if" part is
easy and will be left to the reader.\qed

\vskip10pt

\begin{rem}
The above proof shows that Theorem \ref{OY2} (and so Theorem
\ref{OY3}) is equivalent to saying that the Omori-Yau maximum
principle holds on every Riemannian manifold $(M,g)$ that carries
a positive proper $C^2$ function $\gamma$ satisfying, outside a
compact set,
$$
|\nabla\gamma|\leq u\circ\gamma\;\;\;\;\textnormal{and}
\;\;\;\;\textnormal{Hess}\,\gamma\leq (u\circ\gamma)g
\;\;(\textnormal{resp.} \;\Delta\gamma\leq u\circ\gamma),
$$
where $u\in C^1(0,+\infty)$ is a positive function with $u'\geq 0$
and $\int_1^{+\infty} u(s)^{-1}ds=+\infty$.

\end{rem}

$$
\begin{array}{lccccccl}
\text{Francisco Fontenele}            &&&&&& & \text{Alexandre Paiva Barreto}\\
\text{Departamento de Geometria}      &&&&&& & \text{Departamento de Matem\'atica}\\
\text{Universidade Federal Fluminense}&&&&&& & \text{Universidade Federal de S\~ao Carlos}\\
\text{Niter\'oi, RJ, Brazil}          &&&&&& & \text{S\~ao Carlos, SP, Brazil}\\
\text{fontenele@mat.uff.br}           &&&&&& & \text{alexandre@dm.ufscar.br}\\
\end{array}
$$

\end{document}